  \def\a{\alpha}
  \def\b{\beta}
  \def\d{\delta}
  \def\z{\zeta}
  \def\g{\gamma}
  
  \def\D{\Delta}
  \def\e{\epsilon}
  
   \def\c{\chi}
   \def\p{\psi}

  \def\ra{\rangle}
  \def\ra{\rightarrow}

  \def\go{\rightarrow}
  \def\pf{{\it Proof. }$\;\;$}
  
  \def\hal{\unskip\nobreak\hfil\penalty50\hskip10pt\hbox{}\nobreak
  \hfill\vrule height 5pt width 6pt depth 1pt\par\vskip 2mm}

    \documentclass[11pt,a4paper]{article}
            \usepackage{amsfonts}
            \usepackage{amsmath}

  \DeclareMathOperator{\Aut}{Aut}
\DeclareMathOperator{\Out}{Out}
   \DeclareMathOperator{\Epi}{Epi}
   \DeclareMathOperator{\Irr}{Irr}

\begin{document}

  \title{Commutator maps, measure preservation,
and $T$-systems}
  \author{Shelly Garion
    and
    Aner Shalev \\
    Institute of Mathematics \\
    Hebrew University \\
    Jerusalem 91904 \\Israel}
  \maketitle

  \begin{abstract}
Let $G$ be a finite simple group. We show that the commutator map
$\a:G \times G \go G$ is almost equidistributed as $|G| \go \infty$.
This somewhat surprising result has many applications. It shows that
for a subset $X \subseteq G$ we have $\a^{-1}(X)/|G|^2 = |X|/|G| +
o(1)$, namely $\a$ is almost measure preserving. From this we deduce
that almost all elements $g \in G$ can be expressed as commutators
$g = [x,y]$ where $x,y$ generate $G$.

This enables us to solve some open problems regarding $T$-systems
and the Product Replacement Algorithm (PRA) graph. We show that the
number of $T$-systems in $G$ with two generators tends to infinity
as $|G| \go \infty$. This settles a conjecture of Guralnick and Pak.
A similar result follows for the number of connected components of
the PRA graph of $G$ with two generators.

Some of our results apply for more general finite groups, and more
general word maps.

Our methods are based on representation theory, combining classical
character theory with recent results on character degrees and values
in finite simple groups. In particular the so called Witten zeta
function $\z^G(s) = \sum_{\c \in \Irr(G)} \c(1)^{-s}$ plays a key
role in the proofs.

\end{abstract}

  \footnote{This article was submitted to the {\it Transactions of the
American Mathematical Society} on 21 February 2007 and accepted on
24 June 2007}
  \footnotetext{
  The second author acknowledges the support of grants from the
Israel Science Foundation and the Bi-National Science Foundation
United-States Israel}

\newtheorem{theorem}{Theorem}
  \newtheorem{thm}{Theorem}[section]
  \newtheorem{prop}[thm]{Proposition}
  \newtheorem{lem}[thm]{Lemma}
  \newtheorem{cor}[thm]{Corollary}

\newpage

  \section{Introduction}

\subsection{Finite groups}

Let $G$ be a finite group. Let $\a = \a_G: G \times G \go G$ be the
commutator map, namely
\[
\a(x,y) = [x,y] = x^{-1}y^{-1}xy.
\]
How equidistributed is this map?

To make the question more precise, define for $g \in G$
\[
N(g) = |{\a}^{-1}(g)|,
\]
the size of the fiber above $g$. When can we show that $N(g)$ is
roughly $|G|$ for almost all $g \in G$?

For general groups this is often far from true. However, we shall
show below that commutator maps on finite simple groups are almost
equidistributed. More generally, we associate with each finite group
$G$ a certain parameter $\e(G)$ related to its representation
degrees, and prove that if $\e(G)$ is small then $\a_G$ is almost
equidistributed.

We need some notation. For a finite group $G$ let $\Irr(G)$ denote
the set of complex irreducible characters of $G$. The numbers $N(g)$
above can be studied using a character-theoretic approach, based on
Frobenius classical formula
\begin{equation}
N(g) = |G| \sum_{\c \in \Irr(G)} { {\c(g)} \over {\c(1)} }.
\end{equation}
Thus the character table of $G$ provides complete information on the
distribution of the commutator map. However, computing or estimating
the right hand side of (1) for infinite families of groups, without
complete information on their character tables, is often a
formidable task. We shall see below how information on character
degrees (but not character values) sometimes suffices to draw strong
conclusions regarding commutator maps.

The character degrees of $G$ are conveniently encoded in the so
called \emph{Witten zeta function} $\z^G$ of $G$, defined by
\[
\z^G(s) = \sum_{\c \in \Irr(G)} \c(1)^{-s},
\]
where $s$ is a real number. This function, which plays an important
role in this paper, was originally defined and studied by Witten
\cite{Wit} for Lie groups. For finite simple groups it was studied
and applied in detail in \cite{LiSh3, LiSh4, LiSh5}.

Let $P=P^G$ be the commutator distribution on $G$, namely
\[
P(g) = N(g)/|G|^2,
\]
and let $U$ be the uniform distribution on $G$ (so $U(g) = 1/|G|$
for all $g \in G$).

Our first result bounds the $L_1$-distance
\[
||P-U||_1 = \sum_{g \in G} |P(g)-U(g)|
\]
between the probability measures above.

\begin{prop} Let $G$ be a finite group. Then we have
\[
||P^G-U^G||_1 \le (\z^G(2)-1)^{1/2}.
\]
\end{prop}

Set
\[
\d(G) =  (\z^G(2)-1)^{1/2}.
\]
We now deduce a general lower bound on the number of commutators in
$G$.

\begin{cor} A finite group $G$ has at least
$(1- \d(G))|G|$ commutators.
\end{cor}

Next we define
\[
\e(G) = ( \z^G(2)-1 )^{1/4}.
\]
We establish equidistribution properties of the commutator map on
groups in terms of the parameter $\e(G)$ defined above.

\begin{thm} Every finite group $G$ has a subset
$S = S_G \subseteq G$ with the following properties:

(i) $|S| \ge (1-\e(G))|G|$;

(ii) $(1-\e(G))|G| \le N(g) \le (1+\e(G))|G|$ for all $g \in S$.

\end{thm}

Of course results 1.1-1.3 above have no content when $\z^G(2) \ge 2$
(since then $\d(G), \e(G) \ge 1$). Since the linear characters of
$G$ contribute $|G:G'|$ to $\z^G(2)$ we see that these results can
only be useful for perfect groups, namely groups for which $G = G'$.

Recall that the \emph{representation growth} of $G$ is the series
$\{ r_n(G) \}$, where $r_n(G)$ is the number of irreducible
representations of $G$ of degree $n$. See \cite{LM, Ja} for
background. We clearly have
\[
\z^G(s) = \sum_{n \ge 1} r_n(G)n^{-s}.
\]
If $G$ is perfect, and has very small representation growth, so that
\[
\sum_{n \ge 2} r_n(G)n^{-2} < 1,
\]
then results 1.1-1.3 can be meaningfully applied. Roughly speaking,
they show that perfect groups with few representations have many
commutators, and that commutator maps on them are almost
equidistributed.

This in turn has some further applications.

\begin{cor} Let $G$ be a finite group.

(i) The commutator map $\a = \a_G$ satisfies
\[
\left| { {|\a^{-1}(Y)|} \over {|G|^2} } - { {|Y|} \over {|G|} }
\right| \leq 3\e(G) \hbox{ for all } Y \subseteq G.
\]

(ii) If $X \subseteq G \times G$ then
\[
{ {|\a (X)|} \over {|G|} } \ge { {|X|} \over {|G|^2} } - 3\e(G).
\]
\end{cor}

Part (i) above shows that, if $\e(G)$ is close to zero, then the
commutator map on $G$ is almost measure preserving.

\subsection{Finite simple groups}

The main context in which the above results can be successfully
applied is that of finite simple groups. Indeed, by Theorem 1.1 of
\cite{LiSh4}, if $G$ is simple, then
\[
\z^G(s) \go 1 \hbox{ as } |G| \go \infty \hbox{ provided } s>1.
\]
It follows that for finite simple groups $G$, $\d(G)$ and $\e(G)$
tend to zero as $|G| \go \infty$.

Applying this we deduce our main equidistribution results for
commutator maps. Theorem 1.3 gives rise to the following.

\begin{thm} Every finite simple group $G$ has a subset
$S = S_G \subseteq G$ with the following properties:

(i) $|S| = |G|(1-o(1))$;

(ii) $N(g) = |G|(1+o(1))$ uniformly for all $g \in S$.

\end{thm}

Here and throughout this paper $o(1)$ denotes a real number
depending only on $G$ which tends to zero as $|G| \go \infty$.

The proofs of Theorem 1.5, and of our next results below, rely on
the Classification of Finite Simple Groups. It clearly suffices to
consider alternating groups $A_n$ and simple groups of Lie type. We
show in Section 4 below that $\e(A_n) = O(n^{-1/2})$, and if $G$ is
of Lie type of rank $r$ over a field with $q$ elements then $\e(G) =
O(q^{-r/4})$. This provides explicit upper bounds on the error term
$o(1)$ in Theorem 1.5.

Note that we cannot require in the theorem above that $S=G$. Indeed,
it is well known (and follows from (1) above) that
\[
N(1) = |G|k(G),
\]
where $k(G)$ is the number of conjugacy classes in $G$. Since $k(G)
\go \infty$ as $|G| \go \infty$ we see that the fiber above $g=1$ is
large and does not satisfy condition (ii).

Theorem 1.5 amounts to saying that, for a finite simple group $G$,
\[
||P^G-U^G||_1 \go 0 \hbox{ as } |G| \go \infty.
\]

We present two proofs of Theorem 1.5. The first is probabilistic,
based on Proposition 1.1 bounding the $L_1$-distance above, which
proves the existence of the required subsets $S$. In the second
proof we construct the subsets $S$ explicitly, which often yields
better lower bounds on their cardinality. For example, in the
constructive proof for $A_n$ we obtain $|S| \ge |A_n|(1-
2/[\sqrt{n}]!)$, which is much better than the
$|A_n|(1-O(n^{-1/2}))$ lower bound given by the probabilistic proof.
See Section 5 for more details.

Our constructive proof of Theorem 1.5 applies some powerful recent
results on character values and degrees (see \cite{LiSh3, LiSh4, MS,
Sh}) to estimate the right hand side of Frobenius formula (1) for
specific (almost all) elements $g \in G$, showing that the main
contribution comes from the trivial character $\c = 1$, and the
accumulative contribution of all non-trivial characters is marginal.

Theorem 1.5 can be seen as a culmination of various results on
commutators in finite simple groups. This topic has a long history.
A conjecture of Ore \cite{O} from 1951, which is still not fully
resolved (see \cite{EG}), states that all elements of a finite
simple group are commutators. Various results were obtained in order
to present group elements as commutators or short products of
commutators, see for instance \cite{Ga, Wi, Go, Sh}.

Theorem 1.5 can be viewed as a probabilistic version of Ore's
conjecture. It obviously implies that almost all elements of a
finite simple group are commutators, a result obtained recently in
\cite{Sh}, Theorem 2.9(i).

However, Theorem 1.5 has more refined consequences, as follows.

\begin{cor} Let $G$ be a finite simple group.

(i) The commutator map $\a$ is almost measure preserving, namely
\[
{ {|\a^{-1}(Y)|} \over {|G|^2} } = { {|Y|} \over {|G|} } + o(1)
\hbox{ for all } Y \subseteq G.
\]

(ii) If $X \subseteq G \times G$ then
\[
{ {|\a (X)|} \over {|G|} } \ge { {|X|} \over {|G|^2} } - o(1).
\]

(iii) In particular, if $X \subseteq G \times G$ satisfies $|X| =
(1-o(1))|G|^2$ then
\[
|\a (X)| = (1-o(1))|G|,
\]
namely almost all elements of $G$ can be represented as commutators
$[x,y]$ where $(x,y) \in X$.
\end{cor}

Indeed, this follows by combining Corollary 1.4 with the fact that
$\e(G) = o(1)$.

As a consequence of (ii) above we see that if $A,B$ are subsets of
$G$ such that $|A| \ge a|G|$ and $|B| \ge b|G|$, then the set
$[A,B]$ of commutators $[x,y]$ where $x \in A, y \in B$ satisfies
\[
|[A,B]| \ge (ab - 3\e(G))|G| = (ab-o(1))|G|.
\]

An old conjecture of Dixon \cite{Di}, which is now a theorem (see
\cite{B}, \cite{KL}, \cite{LiSh1}), states that almost all pairs of
elements of a finite simple group are generating pairs. Applying
Corollary 1.6 for the set $X$ of generating pairs of $G$ we obtain

\begin{thm} Let $G$ be a finite simple group, and let $g \in G$
be randomly chosen. Then the probability that $g$ can be represented
as a commutator $g = [x,y]$ where $x,y$ generate $G$ tends to $1$ as
$|G| \go \infty$.
\end{thm}

While Ore's conjecture has been established for alternating groups
and for simple groups of Lie type over large fields, Theorem 1.7 is
new for all types of finite simple groups.

Combining a character-theoretic approach with probabilistic
arguments enables us to obtain additional equidistribution results
(see Theorems 7.1 and 7.4 below). We show that maps on finite simple
groups induced by the word $x^2y^2$, or by longer commutators in any
arrangement of brackets, are almost equidistributed.

\subsection{$T$-systems and the Product Replacement Algorithm}

A main motivation behind Theorem 1.7 comes from the study of
transitivity systems, also known as $T$-systems.

Let $G$ be a finite group and let $d(G)$ be the minimal number of
generators of $G$. For $k \ge 1$ let $F_k$ denote the free group on
$k$ generators. For any $k \geq d(G)$, let
\[
V_k(G)=\{(g_1,\ldots,g_k) \in G^k: \langle g_1,\ldots,g_k\rangle =
G\}
\]
be the set of all generating $k$-tuples of $G$. One can identify
$V_k(G)$ with the set of epimorphisms $\Epi(F_k \ra G)$. The group
$\Aut(F_k) \times \Aut(G)$ acts on $V_k(G)$ by $(\tau, \sigma): \phi
\ra \sigma \circ \phi \circ \tau^{-1}$, where $\tau \in \Aut(F_k),
\sigma \in \Aut(G)$ and $\phi \in \Epi(F_k \ra G)$. The orbits of
this action are called \emph{systems of transitivity}, and also
\emph{$T$-systems} or \emph{$T_k$-systems}, when we specify the
value of $k$. They were introduced by B.H. Neumann and H. Neuman in
\cite{NN} in the context of presentations of finite groups and
studied further in \cite{Du1,Du2,E2,Gi,GP,Ne,P}.

It is well-known that $\Aut(F_k)$ is generated by the following
moves, called the \emph{Nielsen moves}, on the $k$ standard
generators $x_1,\dots,x_k$ of $F_k$, viewed as elements in
$\Aut(F_k)$.
\[
\begin{split}
    R_{i,j}&:  (x_1,\ldots,x_i,\ldots,x_k) \rightarrow
    (x_1,\ldots,x_i\cdot x_j,\ldots,x_k), \\
    L_{i,j}&:  (x_1,\ldots,x_i,\ldots,x_k) \rightarrow
    (x_1,\ldots,x_j \cdot x_i,\ldots,x_k), \\
    P_{i,j}&:  (x_1,\ldots,x_i,\ldots,x_j,\ldots,x_k) \rightarrow
    (x_1,\ldots,x_j,\ldots,x_i,\ldots,x_k),\\
    I_{i}&:  (x_1,\ldots,x_i,\ldots,x_k) \rightarrow
    (x_1,\ldots,x_i^{-1},\ldots,x_k), \\
    \text{for }& 1 \leq i \neq j \leq k.
\end{split}
\]

These moves define the following graph: its vertices are $V_k(G)$
and its edges correspond to the Nielsen moves and to the
automorphisms of $G$. The connected components of this graph are
exactly the $T_k$-systems.

In recent years there is renewed interest in $T$-systems due to
exciting applications to the Product Replacement Algorithm.

The \emph{Product Replacement Algorithm (PRA)} is a practical
algorithm to construct random elements of a finite group. The
algorithm was introduced and analyzed in \cite{CLMNO}, where the
authors proved that it produces asymptotically uniformly distributed
elements. As the success of the algorithm became widely
acknowledged, it was included as a standard routine in the two major
algebra packages GAP and MAGMA. Since then the algorithm was widely
investigated (see \cite{BP,GaP,LP,P}).

The product replacement algorithm is defined as follows
\cite{CLMNO,P}. Given a generating $k$-tuple $(g_1,\ldots,g_k) \in
V_k(G)$, a \emph{move} to another such tuple is defined by first
drawing uniformly a pair $(i,j)$ with $1 \leq i \neq j \leq k$ and
then applying one of the following four operations with equal
probability:
\[
\begin{split}
    R_{i,j}^{\pm}:  (g_1,\ldots,g_i,\ldots,g_k) \rightarrow
    (g_1,\ldots,g_i\cdot g_j^{\pm 1},\ldots,g_k) \\
    L_{i,j}^{\pm}:  (g_1,\ldots,g_i,\ldots,g_k) \rightarrow
    (g_1,\ldots,g_j^{\pm 1} \cdot g_i,\ldots,g_k).
\end{split}
\]

To produce a random element in $G$, start with some generating
$k$-tuple, apply the above moves several times, and finally return a
random element of the generating $k$-tuple that was reached.

The moves in the PRA can be conveniently encoded by the \emph{PRA
graph} $\Gamma_k(G)$ whose vertices are the tuples $V_k(G)$, with
edges corresponding to the moves $R_{i,j}^{\pm}, L_{i,j}^{\pm}$. The
PRA corresponds to a random walk on this graph. However, it is
usually more convenient to look at the following extended graph. The
\emph{extended PRA graph} $\tilde \Gamma_k(G)$ is a graph on
$V_k(G)$ corresponding to the \emph{Nielsen moves}, $R_{i,j}^{\pm},
L_{i,j}^{\pm}$ and $P_{i,j}, I_{i}$, for $1 \leq i \neq j \leq k$.

It is clear from the definitions that the number of $T_k$-systems is
less or equal to the number of connected components in $\tilde
\Gamma_k(G)$, denoted by $\tilde \chi_k(G)$. In addition, if
$\chi_k(G)$ denotes the number of connected components in
$\Gamma_k(G)$, then $\tilde \chi_k(G) \leq \chi_k(G) \leq 2\tilde
\chi_k(G) $. Moreover, if $k\geq d(G)+1$ then $\Gamma_k(G)$ is
connected if and only if $\tilde \Gamma_k(G)$ is connected, and if
$k\geq 2d(G)$, then both $\Gamma_k(G)$ and $\tilde \Gamma_k(G)$ are
connected if and only if $G$ has only one $T_k$-system (see
\cite{P}).

An interesting question is to estimate the number of $T_k$-systems
of $G$ as a function of $k$ and $G$. It is well known that a
$k$-generated abelian group has only one $T_k$-system. A first
example of a nilpotent group $G$ with more than one $T_k$-system,
where $k=d(G)$, was given in \cite{Ne}, as an answer to a question
of Gasch\"utz. Later, Dunwoody \cite{Du1} proved that the number of
$T_k$-systems of certain groups is in fact not bounded, i.e. for
every $k, N$ and $p$, one can find a $p$-group $G$ with $d(G)=k$
such that the number of $T_k$-systems of $G$ is at least $N$.

Particular attention was given to $T$-systems in finite simple
groups $G$. Here $d(G) = 2$ and a conjecture attributed to Weigold
states that for $k \ge 3$ the number of $T_k$-systems of $G$ is 1.
This conjecture was proven for very few families of simple groups
(see Gilman \cite{Gi} and Evans \cite{E}). However, the case $k=2$
seems to be different. It has been shown that the number of
$T_2$-systems in $G = L_2(p) = PSL_2(p)$ tends to infinity as $|G|
\go \infty$ (see Evans \cite{E} and Guralnick and Pak \cite{GP}). A
similar result for  $G = A_n$ was proved by Pak (see \cite{P}). In
\cite{GP} Guralnick and Pak suggest that this might be true for all
finite simple groups, and remark that different methods will have to
be established in order to confirm this.

In this paper we confirm this conjecture. Moreover, we obtain
explicit lower bounds on the number of $T_2$-systems for all finite
simple groups.

For a (possibly twisted) Lie type $L$, not $^2\!B_2,\,^2\!G_2$ or
$^2\!F_4$, define the rank $r = r(L)$ to be the untwisted Lie rank
of $L$ (that is, the rank of the ambient simple algebraic group);
and for $L$ of type $^2\!B_2,\,^2\!G_2$ or $^2\!F_4$, define $r(L) =
1,1,2$ respectively. Let $G_r(q)$ denote a finite simple group of
Lie type of rank $r$ over a field with $q$ elements (in the unitary
case the natural module for $G_r(q)$ is over the field with $q^2$
elements).

\begin{thm} Let $G$ be a finite simple group.
Then the number of $T_2$-systems in $G$ tends to infinity as $|G|
\go \infty$.

Moreover, this number is at least $aq^r r^{-1} (\log{q})^{-2}$ when
$G = G_r(q)$, and at least $n^{(1/2 - \e)\log{n}}$ when $G = A_n$.
\end{thm}

Here $a$ is an absolute positive constant, and $\e>0$ is arbitrary
provided $n$ is large enough (namely $n \ge f(\e)$).

The second, quantitative, assertion of Theorem 1.8, is new even for
$L_2(p)$ and $A_n$, and answers a question from \cite{GP}. The lower
bound for the number of $T_2$-systems in $A_n$ which follows from
the argument in \cite{P} is about $n/(2\log{n})$. No explicit lower
bound was known for $L_2(p)$. In fact the detailed lower bound for
$G_r(q)$ which stems from our proof of Theorem 1.8 is somewhat
better than the one stated, and has the form $(1/8-o(1))p$ when $G =
L_2(p)$.

Theorem 1.8 immediately provides information on the Product
Replacement Algorithm graph in case of two generators.

\begin{cor} Let $G$ be a finite simple group. Then the number of
connected components of the PRA graph $\Gamma_2(G)$ tends to
infinity as $|G| \go \infty$.
\end{cor}

Of course the number of components above is bounded below by the
number of $T_2$-systems, hence the lower bounds of Theorem 1.8
apply. In particular $\Gamma_2(A_n)$ has at least $n^{(1/2 -
\e)\log{n}}$ connected components. For groups of Lie type our method
yields a somewhat better bound, showing that $\Gamma_2(G_r(q))$ has
at least $aq^r (\log{q})^{-1}$ connected components (see Section 6).

The main tool in our proof of Theorem 1.8 is Theorem 1.7, which is
in turn a by-product of our equidistribution theorem (1.5).

\subsection{Notation and layout}

Our notation is rather standard. For a finite group $G$ let
$\Irr(G)$ denote the set of irreducible complex characters of $G$.
Let $k(G)$ denote the number of conjugacy classes in $G$. For $g \in
G$ we let $g^G$ be the conjugacy class of $g$ in $G$, and let $|g|$
be the order of $g$. We denote by $G_r(q)$ a finite simple group of
rank $r$ over a field with $q$ elements. An element $g \in G_r(q)$
is called regular semisimple if it is semisimple and its centralizer
in the ambient algebraic group has minimal dimension $r$. To a
function $f:X \go Y$ between finite sets we associate a probability
distribution $P_f$ on $Y$ given by
\[
P_f(y) = \frac{|f^{-1}(y)|}{|X|} \hbox{ where } y \in Y.
\]
By a word $w = w(x_1, \dots , x_m)$ we mean an element of the free
group $F_m$ on $x_1, \dots , x_m$. Given a group $G$ the word $w$
defines a function $G^m \go G$, obtained by substitution, which we
denote by $\a_w$. The associated probability distribution $P_{\a_w}$
on $G$ will be also denoted by $P_w(G)$. Additional notation will be
introduced when needed.

Some words on the layout of this paper. In Section 2 we use Fourier
techniques to prove Proposition 1.1. In Section 3 we prove results
1.2-1.4. Section 4 presents recent delicate results on simple groups
which are required in order to prove Theorem 1.5 constructively. The
constructive proof of this theorem is then carried out in Section 5.
Section 6 is devoted to the various applications, focusing on
$T$-systems; this is where results 1.7-1.9 are proved. In Section 7
we obtain additional equidistribution results for longer
commutators.

\emph{Acknowledgments.} We are grateful to Benjy Weiss for useful
remarks. This paper is part of the first author's Ph.D. thesis done
under the supervision of Alex Lubotzky.

\bigskip
\bigskip

\section{Bounding the $L_1$-distance}

The main purpose of this section is to prove Proposition 1.1. We
shall start with a more general discussion, from which this result
will be deduced.

Let $G$ be a finite group, and let $P$ be a probability distribution
on $G$ which is a class function. For example, we may have $P = P_w$
where $w$ is a word.

Consider the non-commutative Fourier expansion
\[
P = |G|^{-1} \sum_{\c \in \Irr(G)} a_{\c} \c,
\]
with suitable (complex) coefficients $a_{\c}$.

Note that $a_1 = 1$ (since $\sum_{g \in G} P(g) = 1$).

\begin{lem} We have
\[
\sum_{g \in G} P(g)^2 = |G|^{-1} \sum_{\c \in \Irr(G)} |a_{\c}|^2.
\]
\end{lem}

\pf We have
\[
\sum_{g \in G} P(g)^2 = \sum_{g \in G} P(g) {\overline{P(g)}} =
|G|^{-2} \sum_{g \in G} \sum_{\c} a_{\c}\c(g) \sum_{\p}
{\overline{a_{\p}}} {\overline{\p(g)}}
\]
\[
= |G|^{-2} \sum_{\c,\p} a_{\c} {\overline{a_{\p}}} \sum_{g \in G}
\c(g){\overline{\p(g)}}.
\]
Using the orthogonality relations this gives
\[
\sum_{g \in G} P(g)^2 = |G|^{-1} \sum_{\c,\p}
a_{\c}{\overline{a_{\p}}} \d_{\c \p} = |G|^{-1} \sum_{\c}
|a_{\c}|^2.
\]
\hal
\medskip

\begin{lem} We have
\[
\sum_{g \in G} (P(g)-|G|^{-1})^2 = |G|^{-1} \sum_{\c \ne 1}
|a_{\c}|^2.
\]
\end{lem}

\pf Using the previous lemma we have
\[
\sum_{g \in G} (P(g)-|G|^{-1})^2 = \sum_{g \in G} P(g)^2 - 2|G|^{-1}
\sum_{g \in G} P(g) + |G|^{-1} = |G|^{-1}  (\sum_{\c} |a_{\c}|^2 -
1).
\]
The result follows since $a_1 = 1$. \hal
\medskip

\begin{lem} We have
\[
||P-U||_1 \le  (\sum_{\c \ne 1} |a_{\c}|^2)^{1/2}.
\]
\end{lem}

\pf By Cauchy-Schwarz inequality we have
\[
(||P-U||_1)^2 = (\sum_{g \in G} |P(g)-|G|^{-1}|)^2 \le |G| \sum_{g
\in G} (P(g)-|G|^{-1})^2.
\]
The result follows using the lemma above. \hal
\medskip

\noindent {\bf Proof of Proposition 1.1}

Let $P$ be the commutator probability distribution on $G$. Then by
(1) we have
\[
a_{\c} = \c(1)^{-1} \hbox{ for all } \c.
\]
Applying the above lemma we obtain
\[
||P-U||_1 \le (\sum_{\c \ne 1} \c(1)^{-2})^{1/2} =
(\z^G(2)-1)^{1/2}.
\]
This completes the proof. \hal
\medskip

Using the methods of this section we can easily show the following.

\begin{prop} Let $G$ be a finite group, and let $P_{x^2y^2}$
be the probability distribution associated with the word map from $G
\times G$ to $G$ induced by $x^2y^2$. Let $R$ denote the set of real
characters of $G$. Then
\[
||P_{x^2y^2}-U||_1 \le (\sum_{\c \in R} \c(1)^{-2} -1)^{1/2} \le
(\z^G(2)-1)^{1/2}.
\]
\end{prop}

\pf Let $\b = \a_{x^2y^2}:G \times G \go G$. We use the well known
formula (see e.g. \cite{LiSh3}, 3.1)
\[
|\b^{-1}(g)| = |G| \sum_{\c \in R} {{\c(g)} \over {\c(1)}}.
\]
This means that
\[
P_{\b} = |G|^{-1} \sum_{\c \in R} \c(1)^{-1} \c.
\]
The result now follows from Lemma 2.3. \hal
\medskip

\bigskip
\bigskip

\section{Applications of the $L_1$ bound}

The purpose of this section is to deduce results 1.2-1.4 from
Proposition 1.1.

\noindent {\bf Proof of Corollary 1.2}

Set $\d = \d(G)$. Then $||P_f - U||_1 \le \d$ by Proposition 1.1.
Let $C$ be the set of commutators in $G$, and let $D = G \setminus
C$. Then
\[
\d \ge \sum_{g \in D} \left| P(g) - \frac{1}{|G|} \right|
  =  |D|/|G|.
\]
Therefore $|D| \le \d |G|$, so $|C| \ge (1-\d)|G|$, as required.
\hal
\medskip

The notion of almost equidistribution is naturally defined in the
general setting of arbitrary functions between finite sets. Let
$X,Y$ be finite sets, and let $\e >0$. We say that a function $f:X
\ra Y$ is \emph{$\e$-equidistributed} if there exists a subset $Y'
\subseteq Y$ with the following properties:

(i) $|Y'| \geq |Y|(1-\e)$;

(ii) $\frac{|X|}{|Y|}(1-\e) \leq |f^{-1}(y)| \leq
\frac{|X|}{|Y|}(1+\e)$ uniformly for all $y \in Y'$.

Recall that $P_f$ is the probability distribution on $Y$ induced by
$f$. Let $U$ be the uniform distribution on $Y$. We show that the
$L_1$-distance
\[
    || P_f - U ||_1 = \sum_{y \in Y} \left|
    P_f(y) - \frac{1}{|Y|} \right|
\]
is small if and only if $f$ is almost equidistributed. Indeed we
have the following easy Lemma.

\begin{lem} With the above notation we have

(i) If $f$ is $\e$-equidistributed then $||P_f-U||_1 \le 4 \e$.

(ii) If $||P_f - U||_1 \le \d$, then $f$ is
$\sqrt{\d}$-equidistributed.

\end{lem}

\pf (i) Assume that $f:X \rightarrow Y$ is $\e$-equidistributed, and
let $Y'$ be as in the definition above. Write $Y = Y' \cup Y''$ as a
disjoint union. Then for any $y \in Y'$, $\left|P_f(y) -
\frac{1}{|Y|}\right| \leq \frac{\e}{|Y|}$, and thus
\[
    \left| \sum_{y \in Y'} P_f(y) - \sum_{y \in Y'} \frac{1}{|Y|} \right|
    \leq \sum_{y \in Y'} \left| P_f(y) - \frac{1}{|Y|} \right| \leq
    |Y'|\frac{\e}{|Y|} \leq \e.
\]
However, $\sum_{y \in Y'} \frac{1}{|Y|} = \frac{|Y'|}{|Y|} \geq
1-\e$, therefore $\sum_{y \in Y'} P_f(y) \geq 1-2\e$, so we deduce
that $\sum_{y \in Y''} P_f(y) \leq 2\e$. Therefore
\begin{align*}
    ||P_f - U||_1 &= \sum_{y \in Y} \left|
    P_f(y) - \frac{1}{|Y|} \right| = \sum_{y \in Y'} \left|
    P_f(y) - \frac{1}{|Y|} \right| + \sum_{y \in Y''} \left|
    P_f(y) - \frac{1}{|Y|} \right| \\
   &\leq \sum_{y \in Y'} \frac{\e}{|Y|} + \sum_{y \in
   Y''}|P_f(y)| + \sum_{y \in Y''}\frac{1}{|Y|} \\
   &\leq |Y'|\frac{\e}{|Y|} + 2\e + |Y''|\frac{1}{|Y|}
   \leq \e + 2\e + \e = 4\epsilon .
\end{align*}

(ii) Assume that $||P_f - U||_1 \le \d$. Define $Y'' = \left\{y \in
Y: \left|
    P_f(y) - \frac{1}{|Y|} \right| > \frac{\sqrt{\d}}{|Y|}
    \right\}$. Then
\begin{align*}
    \d \geq \sum_{y \in Y} \left| P_f(y) - \frac{1}{|Y|} \right|
    \geq \sum_{y \in Y''} \left| P_f(y) - \frac{1}{|Y|} \right|
    > |Y''| \frac{\sqrt\d}{|Y|}.
\end{align*}
Therefore, $|Y''| < \sqrt{\d}|Y|$. Take $Y'=Y \setminus Y''$. Then
$|Y'| \geq |Y|(1-\sqrt{\d})$ and any $y \in Y'$ satisfies $\left|
P_f(y) - \frac{1}{|Y|} \right| \leq \frac{\sqrt{\d}}{|Y|}$. Thus,
$f$ is $\sqrt{\d}$-equidistributed. \hal
\medskip

\noindent {\bf Proof of Theorem 1.3}

This follows by combining Proposition 1.1 and part (ii) of Lemma
3.1. \hal
\medskip

The next result concerns measure preservation.

\begin{prop}
Let $f: X \ra Y$ be $\e$-equidistributed.

(i) If $Y_0 \subseteq Y$ then
\[
    \left|\frac{f^{-1}(Y_0)}{|X|} - \frac{|Y_0|}{|Y|} \right| \leq
    3\e.
\]

(ii) If $X_0 \subseteq X$ then
\[
    \frac{|f(X_0)|}{|X|} \geq \frac{|X_0|}{|X|} - 3\e.
\]
\end{prop}

\pf Assume that $f: X \ra Y$ is $\e$-equidistributed and let $Y'$ be
as in the definition above. Let $X' = f^{-1}(Y') \subseteq X$ be the
inverse image of $Y'$. Then by part (i) and the lower bound in (ii)
of the definition,
\[
    |X'| = \sum_{y \in Y'} |f^{-1}(y)| \geq |Y'|
    \frac{|X|}{|Y|}(1-\e) \geq |X|(1-\e)^2 \geq |X|(1-2\e).
\]
We conclude that
\[
    |f^{-1}(Y \setminus Y')| = |X| - |X'| \leq 2\e |X|.
\]

Now let $Y_0 \subseteq Y$. Then
\[
    |f^{-1}(Y_0)| \leq |f^{-1}(Y_0 \cap Y')| + |f^{-1}(Y \setminus
    Y')| \leq \sum_{y \in Y_0 \cap Y'}|f^{-1}(y)| + 2\e |X|.
\]
Using the upper bound in part (ii) of the definition, we see that
\[
    |f^{-1}(Y_0)| \leq |Y_0|\frac{|X|}{|Y|}(1+\e) + 2\e |X|.
\]
Therefore
\[
    \frac{|f^{-1}(Y_0)|}{|X|} \leq \frac{|Y_0|}{|Y|}(1+\e) + 2\e \leq
    \frac{|Y_0|}{|Y|} + 3\e.
\]

On the other hand we have
\begin{align*}
    |f^{-1}(Y_0)| &\geq |f^{-1}(Y_0 \cap Y')| = \sum_{y \in Y_0 \cap
    Y'} |f^{-1}(y)| \geq |Y \cap Y_0| \frac{|X|}{|Y|}(1-\e) \\ &\geq (|Y_0| -
    \e|Y|)\frac{|X|}{|Y|}(1-\e).
\end{align*}
This yields
\[
    \frac{|f^{-1}(Y_0)|}{|X|} \geq \left(\frac{|Y_0|}{|Y|} - \e
    \right)(1-\e) \geq \frac{|Y_0|}{|Y|} - 2\e.
\]

This completes the proof of part (i) of the proposition.

Part (ii) now follows easily. Indeed, given $X_0 \subseteq X$,
define $Y_0 = f(X_0)$. As above, we have
\[
    \frac{|f^{-1}(Y_0)|}{|X|} \leq \frac{|Y_0|}{|Y|} + 3\e =
    \frac{|f(X_0)|}{|Y|} + 3\e.
\]
Since $X_0 \subseteq f^{-1}(Y_0)$ we obtain
\[
    \frac{|f(X_0)|}{|Y|} \geq \frac{|f^{-1}(Y_0)|}{|X|} -3\e \geq \frac{|X_0|}{|X|} - 3\e.
\]
\hal
\medskip

\noindent {\bf{Proof of Corollary 1.4}}

By Theorem 1.3 the commutator map $\a: G \times G \ra G$ is
$\e(G)$-equidistributed. The corollary now follows immediately from
Proposition 3.2.

\hal
\medskip

\bigskip
\bigskip

\section{Simple groups: character theoretic tools}

In this section we introduce the main concepts and tools which are
needed for our two proofs of Theorem 1.5.

The first result summarizes some of the properties of the Witten
zeta function $\z^G$ defined in the Introduction. Our first proof of
Theorem 1.5 follows by combining part (i) below with Theorem 1.3.

\begin{thm} Let $G$ be a finite simple group.

(i) If $s>1$ then $\z^G(s) \go 1$ as $|G| \go \infty$.

(ii) If $s>2/3$ and $G \ne L_2(q)$ then $\z^G(s) \go 1$ as $|G| \go
\infty$.

(iii) If $s>0$ and $G = A_n$ then $\z^G(s) \go 1$ as $|G| \go
\infty$. Moreover, $\z^G(s) = 1 + O(n^{-s})$.

(iv) If $G = G_r(q)$ then $\z^G(2) = 1 + O(q^{-r})$.
\end{thm}

Indeed, parts (i) and (ii) are Theorem 1.1 of \cite{LiSh4}. Part
(iii) is Corollary 2.7 of \cite{LiSh3}.

To prove Part (iv) we note that
\[
\z^G(2) - 1 \le k(G)h(G)^{-2},
\]
where $h(G)$ is the minimal degree of a non-trivial character of
$G$. Suppose $G = G_r(q)$. Then
\[
h(G) \ge c_1 q^r,
\]
by Landazuri-Seitz \cite{LaSe}, and
\[
k(G) \le c_2 q^r
\]
by Fulman-Guralnick \cite{FG}, where $c_1, c_2 > 0$ are absolute
constants. Part (iv) now follows from the above inequalities.

The next result we need deals with regular semisimple elements in
groups of Lie type.

\begin{thm} Let $G = G_r(q)$ and let $S$ be the set of regular
semisimple elements in $G$.

(i) There is an absolute constant $a_1$ such that
\[
|S| \ge |G|(1-a_1 q^{-1});
\]

(ii) There is a number $f(r)$ depending only on $r$ such that, if $g
\in S$, then
\[
|\c(g)| \le f(r)
\]
for all $\c \in \Irr(G)$.

\end{thm}

Indeed, part (i) is a result of Guralnick and L\"ubeck \cite{GL},
while part (ii) is Lemma 4.4 of \cite{Sh}.

The next result we use deals with elements $g \in G_r(q)$ whose
centralizer is not very large.

\begin{thm} Let $G = G_r(q)$, and fix $\e$ satisfying
$ 0 < \e < 2$. Let
\[
S(\e) = \{ g \in G: |C_G(g)| \le q^{(3-\e)r} \}.
\]
Then we have

(i) $|S(\e)| \ge |G|(1 - a_2 q^{-(2-\e)r})$, where $a_2$ is an
absolute constant;

(ii) There is a number $r_1(\e)$ such that, if $r \ge r_1(\e)$, and
$g \in S(\e)$, then
\[
\sum_{1 \ne \c \in \Irr(G)} {{|\c(g)|} \over {\c(1)}} \go 0 \hbox{
as } |G| \go \infty.
\]
\end{thm}

Indeed part (i) follows from Corollary 5.4 of \cite{Sh}, while part
(ii) is Proposition 4.7 of \cite{Sh}.

We conclude this section quoting a useful result of M\"uller and
Schlage-Puchta on character values for symmetric groups. See part
(i) of Theorem A in \cite{MS}.

\begin{thm} Let $g \in S_n$ be a permutation with $f$ fixed
points. Define
\[
\d = ( (1-1/\log{n})^{-1} { {12 \log{n}} \over \log{(n/f)} } + 18
)^{-1}.
\]
Then we have
\[
|\c(g)| \le \c(1)^{1-\d}
\]
for all $\c \in \Irr(S_n)$.
\end{thm}

\medskip
\noindent {\bf Remark}

By Theorem 4.1 we have $\z^G(2)-1 = O(n^{-2})$ where $G = A_n$, and
$\z^G(2)-1 = O(q^{-r})$ where $G = G_r(q)$. This shows that
\[
\e(A_n) = O(n^{-1/2}) \hbox{ and } \e(G_r(q)) = O(q^{-r/4}).
\]
\noindent

\medskip

\bigskip
\bigskip

\section{Theorem 1.5: constructive proof}

In this section we prove Theorem 1.5 in a constructive manner,
providing the subsets $S \subset G$ with the required properties. We
need some notation.

For $g \in G$ let
\[
\D(g) = \sum_{1 \ne \c \in \Irr(G)} { {\c(g)} \over {\c(1)} },
\]
and
\[
E(g) = \sum_{1 \ne \c \in \Irr(G)} { {|\c(g)|} \over {\c(1)} }.
\]

Then $N(g) = |G|(1+\D(g))$ and $|\D(g)| \le E(g)$.

To prove Theorem 1.5 constructively it therefore suffices to find
subsets $S \subseteq G$ consisting of almost all elements of $G$
such that $E(g) = o(1)$ for all $g \in S$. We will see below that
this strategy works out for some (but not all) finite simple groups
$G$. For the remaining groups we will show that $|\D(g)| = o(1)$,
which also suffices.

Our study of $E(g)$ and $\D(g)$ is based on recent detailed results
on character degrees and values quoted in Section 4.

{\bf Case 1.} $G = L_2(q)$.

The character table of $L_2(q)$ is well known (see for instance
\cite{Do}), so $\D(g)$ and $N(g)$ can be computed for all $g \in G$
using (1). We summarize the result below.



\begin{prop}
Let $q=p^n$ be a prime power, then in the group $L_2(q)$,
\begin{enumerate}
\item
If $a$ is an element of order $\frac{q-1}{2}$ (when $q$ is odd) or
order $q-1$ (when $q$ is even) then $\Delta(a^l) = \frac{1}{q} +
\frac{\alpha(q,l)}{q+1}$ where
\[
\alpha(q,l) = \begin{cases} 2 &\text{ if } q \equiv 1\pmod 4, \quad
$l$ \text{ even } \\ -4 &\text{ if } q \equiv 1\pmod 4, \quad $l$
\text{ odd } \\ -1 &\text{ if } q \equiv 3\pmod 4 \text{ or } $q$
\text{ is even}
\end{cases}
\]
\item
If $b$ is an element of order $\frac{q+1}{2}$ (when $q$ is odd) or
order $q+1$ (when $q$ is even) then $\Delta(b^m) = - \frac{1}{q} +
\frac{\beta(q,m)}{q-1}$ where
\[
\beta(q,m)= \begin{cases} 1 &\text{ if } q \equiv 1\pmod 4 \text{ or
} $q$ \text{ is even} \\ -2 &\text{ if } q \equiv 3\pmod 4, \quad
$m$ \text{ even } \\ 4 &\text{ if } q \equiv 3\pmod 4, \quad $m$
\text{ odd }
\end{cases}
\]
\item
If $c$ is an element of order $p$ then
\[
\Delta(c) = \begin{cases} \frac{1}{2(q+1)} & \text{ if } q \equiv 1
\pmod 4 \\  -\frac{1}{q+1}
-\frac{3}{2(q-1)} & \text{ if } q \equiv 3 \pmod 4 \\
-\frac{3}{2(q+1)}-\frac{1}{2(q-1)} & \text{ if } q \text{ is even }
\end{cases}
\]
\end{enumerate}
\end{prop}


Letting $S$ be all non-identity elements of $G$ we find by
Proposition 5.1 that for $g \in S$,
\[
|\D(g)| \le 5/q = o(1).
\]

This implies the required conclusion.

\medskip

Let $\e, r_1(\e)$ be as in Theorem 4.3, and set $b = r_1(1)$.

\medskip

{\bf Case 2.} $G = G_r(q)$, where $r < b$, and $G \ne L_2(q)$.

Let $S$ be the set of regular semisimple elements of $G$. Since the
rank $r$ is bounded, we have $q \go \infty$ as $|G| \go \infty$.

By Theorem 4.2(i) above we have
\[
|S| \ge |G|(1-a_1q^{-1}) = |G|(1-o(1)).
\]

By part (ii) of Theorem 4.2 there is a number $f(r)$ depending only
on $r$ such that
\[
|\c(g)| \le f(r)
\]
for all $\c \in \Irr G$ and $g \in S$. This yields
\[
E(g) \le \sum_{\c \ne 1} f(r)\c(1)^{-1} = f(r)(\z^G(1)-1).
\]

Since $G \ne L_2(q)$ we have $\z^G(1) \go 1$ as $|G| \go \infty$
(see part (ii) of Theorem 4.1). This yields $E(g) = o(1)$ uniformly
for all $g \in S$, proving the result in this case.

\medskip

{\bf Case 3.} $G = G_r(q)$ where $r \ge b$.

We adopt the notation of Theorem 4.3 and apply this result with $\e
= 1$. Part (i) then yields
\[
|S(1)| \ge |G|(1-a_2 q^{-r}).
\]
Now let $g \in S(1)$. Since $r \ge b = r_1(1)$, part (ii) of Theorem
4.3 shows that $E(g) = o(1)$, and again the conclusion follows with
$S = S(1)$.

\medskip

{\bf Case 4.} $G = A_n$.

Let $S$ be the set of permutations in $A_n$ with at most $\sqrt{n}$
fixed points.

It is easy to see that the probability that a permutation $g \in
A_n$ has at least $f$ fixed points is at most $2/f!$. This implies
that
\[
|S| \ge |A_n|(1- 2/[\sqrt{n}]!) = |A_n|(1-o(1)).
\]

Now set
\[
\d = { 1 \over {43} },
\]
and let $g \in S$. Using Theorem 4.4 above we see that, for $n$
large, we have
\begin{equation}
|\c(g)| \le \c(1)^{1-\d}
\end{equation}
for all $\c \in \Irr(S_n)$.

For each irreducible character $\c$  of $S_n$, either $\c \downarrow
A_n$ is irreducible, or $\c \downarrow A_n = \c_1 + \c_2$, a sum of
two irreducible characters of degree $\c(1)/2$. All irreducible
characters of $A_n$ occur in this way.

In the latter case, note that
\[
{ {\c_1(g)} \over {\c_1(1)} } + { {\c_2(g)} \over {\c_2(1)} } = 2 {
{\c(g)} \over {\c(1)} }.
\]
This implies that
\[
\left| \sum_{1 \ne \c \in \Irr(A_n)} { {\c(g)} \over {\c(1)} }
\right| \le 2 \sum_{\c \in \Irr(S_n), \c(1)>1} { {|\c(g)|} \over
{\c(1)} } \le 2 \sum_{\c \in \Irr(S_n), \c(1)>1} \c(1)^{-\d},
\]
where the last inequality follows from (2). We conclude that, in
$A_n$ we have
\[
|\D(g)| \le 2(\z^{S_n}(\d)-2).
\]
By Theorem 1.1 in \cite{LiSh3}, $\z^{S_n}(\d) = 2 + O(n^{-\d})$.
This yields
\[
|\D(g)| = O(n^{-1/43}) \go 0 \hbox{ as } n \go \infty.
\]

This completes the case $G = A_n$ and the constructive proof of
Theorem 1.5. \hal
\medskip

\noindent {\bf Remark}

The proof above yields explicit lower bounds on $|S|$. In some cases
they may be further improved. Indeed, fix any $\e>0$. Then for $G =
G_r(q)$ and $r \ge r_1(\e)$ we can use Theorem 4.3 to construct $S$
satisfying
\[
|S| \ge |G|(1 - a_2q^{-(2-\e)r}).
\]
For $G = A_n$ we may take $S$ as all even permutations with at most
$n^{1-\e}$ fixed points, obtaining
\[
|S| \ge |A_n|(1- 2/[n^{1-\e}]!) \ge |A_n|(1 - n^{-n^{\gamma}}),
\]
where $\gamma < 1$ is arbitrarily close to 1.

\bigskip
\bigskip

\section{Applications to $T$-systems}

In this section we focus on the various applications and prove
results 1.7-1.9.

\noindent {\bf Proof of Theorem 1.7}

This is an immediate consequence of Corollary 1.6(iii), in view of
the fact that almost all pairs of finite simple groups are
generating pairs (see \cite{LiSh1} and the references therein).

\hal
\medskip

We now turn to the main applications, involving $T_2$-systems.

\noindent {\bf{Proof of Theorem 1.8}}

It suffices to prove the second (quantitative) assertion in the
theorem.

Higman's Lemma states that for $k=d(G)=2$, the union of the
conjugacy classes under $\Aut(G)$ of the commutators $[g_1,g_2]$ and
$[g_2,g_1]$ is an invariant of the $T_2$-system of $(g_1,g_2)$ (see
\cite{GP,Ne,P}). This lemma becomes useful when dealing with
$T_2$-systems of finite simple groups $G$.

By Theorem 1.7 there is a subset $S \subset G$ such that $S =
|G|(1-o(1))$, and every $g \in S$ can be written as $g = [g_1,g_2]$
where $(g_1,g_2) \in V_2(G)$ (namely $g_1,g_2$ generate $G$, thus
giving rise to a $T_2$-system).

Let $k(S)$ denote the number of distinct unions $C \cup C^{-1}$
where $C$ is an $\Aut(G)$-conjugacy class of an element of $S$. Then
the number of $T_2$-systems in $G$ is at least $k(S)$.

Suppose first $G = A_n$. Then $\Aut(G) = S_n$ and $C = C^{-1}$ for
any $S_n$-class $C$. Given $\e > 0$ fix $\d > 0$ such that $\d <
\e$. We may assume $n$ is sufficiently large (given $\e$). By a
result of Erd{\H o}s and Tur{\' a}n (see \cite{ET}), the order of
almost all permutations in $S_n$ is at least $n^{(1/2 -
\d)\log{n}}$. Intersecting $S$ with the set of permutations with
this order restriction we may therefore assume that each $g \in S$
has order at least $n^{(1/2 - \d)\log{n}}$, while we still have $|S|
= |A_n|(1-o(1)) \ge n!/3$.

Now, if $g \in S$, then we have
\[
|C_{S_n}(g)| \ge |g| \ge n^{(1/2 - \d)\log{n}}.
\]
Thus
\[
|g^{S_n}| \le n! \cdot n^{-(1/2 - \d)\log{n}}.
\]
Since the union of the $S_n$-classes of $g \in S$ covers $S$ we have
\[
k(S) \cdot n! \cdot n^{-(1/2 - \d)\log{n}} \ge |S| \ge n!/3.
\]
This yields
\[
k(S) \ge  n^{(1/2 - \d)\log{n}}/3 \ge  n^{(1/2 - \e)\log{n}},
\]
proving the result for $A_n$.

Now suppose $G = G_r(q)$.

Define
\[
c(G) = \min \{ |C_G(g)|: g \in G \}.
\]
Let $\Out(G) = \Aut(G)/G$ and let $g \in G$. Then
\[
|g^{\Aut(G)}| = |\Aut(G)|/|C_{\Aut(G)}(g)| \le |\Aut(G)|/c(G).
\]
In particular the size of each union $C \cup C^{-1}$ where $C =
g^{\Aut(G)}$ and $g \in S$ is at most $2 |\Aut(G)|/c(G)$. Since
there are $k(S)$ such unions, and their union covers $S$, we see
that
\[
k(S) \cdot 2 |\Aut(G)|/c(G) \ge |S| \ge |G|(1-o(1)).
\]
This yields
\[
k(S) \ge (1/2-o(1)) c(G)/|\Out(G)|.
\]

By results of Fulman and Guralnick \cite{FG} there is an absolute
constant $c_1>0$ such that
\[
c(G_r(q)) \ge c_1q^r/\log{q}.
\]
The structure of $\Out(G)$ is known, and assuming $q = p^f$ ($p$
prime) we have
\[
|\Out(G)| \le c_2 rf \le c_3 r \log{q}
\]
for some (small) constant $c_3$. The above inequalities yield
\[
k(S) \ge c_4 q^{r} r^{-1} (\log{q})^{-2}.
\]

This completes the proof of Theorem 1.8. \hal
\medskip

\noindent {\bf Remarks}

1. Note that, by \cite{ET}, almost all permutations in $S_n$ have at
most $(1+ \d)\log{n}$ cycles. These permutations split into at most
$n^{(1+ \e)\log{n}}$ conjugacy classes (since this number bounds the
number of partitions of $n$ into at most $(1+ \d)\log{n}$ parts,
where $\e > \d$). Thus, although $k(S_n) = p(n) > c^{\sqrt{n}}$, a
union of just $n^{(1+ \e)\log{n}}$ conjugacy classes covers almost
all of $S_n$. This shows that our lower bound on $k(S)$ in the proof
above is essentially best possible.

2. If $G = G_r(q)$ the proof above shows that the number of
$T_2$-systems in $G$ is at least  $(1/2-o(1)) c(G)/|\Out(G)|$. This
produces specific bounds which are slightly better than the general
one stated in the theorem. For example, it follows that the number
of $T_2$-systems in $L_2(p)$ is at least $(1/8-o(1))p$.
\medskip

We conclude this section by briefly discussing the PRA graph
$\Gamma_2(G)$. An elementary calculation shows that the conjugacy
class in $G$ of the commutator $[g_1,g_2]$ is an invariant of the
connected component of $(g_1,g_2)$ in $\Gamma_2(G)$. Arguments
similar to the proof of Theorem 1.8 show that the number of
components of this graph is at least $(1-o(1))c(G)$. This shows that
$\Gamma_2(G_r(q))$ has at least $a q^r /\log{q}$ connected
components, giving a better lower bound in Corollary 1.9. The lower
bound for $A_n$ remains the same.

\bigskip
\bigskip

\section{Equidistribution revisited}

Any word $w(x_1, \ldots , x_m) \in F_m$ gives rise to a word map
$\a_w : G^m \go G$. Word maps on algebraic groups and on finite
simple groups have been the subject of active investigations in
recent years, see \cite{Bo}, \cite{LiSh2}, \cite{La}, \cite{Sh} and
\cite{LaSh}.

It is interesting to find out which words $w$ have the remarkable
property of the commutator map, namely that $\a_w$ is almost
equidistributed (namely $o(1)$-equidistributed) on all finite simple
groups $G$.

In general this is highly unexpected. For example, power words $w =
x_1^k$ ($k \ge 2$), or words which are proper powers $w = w_1^k$,
have image much smaller than $(1-o(1))|G|$ for infinite families of
finite simple groups $G$, hence their associated maps are not almost
equidistributed.

However, we do obtain positive result for some more words.

\begin{thm} Let $G$ be a finite simple group, and let
$\beta: G \times G \go G$ be the map given by $\b(x,y) = x^2y^2$.
Then there is a subset $S \subset G$ with $|S| = (1-o(1))|G|$ such
that $|\beta^{-1}(g)| = (1+o(1))|G|$ for all $g \in S$.
\end{thm}

\pf By Proposition 2.4 $||P_{\b}-U||_1 \le \d(G)$. Lemma 3.1(ii) now
shows that $\b$ is $\e(G)$-distributed. Finally, by Theorem 4.1(i),
$\e(G) = o(1)$. The result follows. \hal
\medskip

To obtain more positive results we need some preparations.

The following two lemmas show that the property of almost
equidistribution behaves well under direct products and
compositions. The proofs use the $L_1$ notation, which is more
natural.

\begin{lem}
Let $X_1,X_2,Y_1,Y_2$ be finite sets and let $\delta_1,\delta_2>0$.
For $i=1,2$ denote by $U_i$ the uniform distribution on $Y_i$ and
assume that $f_i:X_i \rightarrow Y_i$ satisfies $||P_{f_i}-U_i||_1
\le \d_i$.

Then the function $f=f_1 \times f_2: X_1 \times X_2 \rightarrow Y_1
\times Y_2$, which is defined by
\[
    f(x_1,x_2) = (f_1(x_1),f_2(x_2)) \text{ for } x_1 \in X_1, x_2 \in
    X_2,
\]
satisfies $||P_f-U||_1 \le \d_1+\d_2$, where $U$ is the uniform
distribution on $Y_1 \times Y_2$.
\end{lem}

\pf We have $P_{f_1 \times f_2}(y_1,y_2) = P_{f_1}(y_1)
P_{f_2}(y_2)$ for any $(y_1,y_2) \in Y_1 \times Y_2$. Thus,
\begin{align*}
    || P_{f_1 \times f_2} - U||_1
    &= \sum_{(y_1,y_2) \in Y_1 \times Y_2}
    \left| P_{f_1}(y_1) P_{f_2}(y_2) - \frac{1}{|Y_1|}\frac{1}{|Y_2|}
    \right| \\
    &\leq \sum_{y_2 \in Y_2} P_{f_2}(y_2) \sum_{y_1 \in Y_1} \left|
    P_{f_1}(y_1) - \frac{1}{|Y_1|} \right| \\ &\quad + \sum_{y_1 \in Y_1}
    \frac{1}{|Y_1|} \sum_{y_2 \in Y_2} \left| P_{f_2}(y_2) - \frac{1}{|Y_2|}
    \right| \\
    &\leq \sum_{y_2 \in Y_2} P_{f_2}(y_2) \delta_1 + \sum_{y_1 \in
    Y_1} \frac{1}{|Y_1|} \delta_2 = \delta_1 + \delta_2.
\end{align*}
\hal
\medskip

\begin{lem}
Let $X,Y,Z$ be finite sets, and let $\delta_1,\delta_2 > 0$. Denote
by $U_Y$ and $U_Z$ the uniform distributions on $Y$ and $Z$
respectively. Assume that $f_1:X \rightarrow Y$ and $f_2:Y
\rightarrow Z$ satisfy $||P_{f_1}-U_Y||_1 \le \d_1$ and
$||P_{f_2}-U_Z||_1 \le \d_2$. Then their composition $f = f_2 \circ
f_1: X \rightarrow Z$ satisfies $||P_{f}-U_Z||_1 \le \d_1+\d_2$.
\end{lem}

\pf Note that for any $z \in Z$,
\[
    P_{f_2 \circ f_1}(z) = \frac{|f_1^{-1} (f_2^{-1}(z))|}{|X|} =
    \sum_{f_2(y)=z} \frac{|f_1^{-1}(y)|}{|X|} =
    \sum_{f_2(y)=z}P_{f_1}(y).
\]
Thus, by the assumption on $f_1$,
\begin{align*}
    \sum_{z \in Z} \left| P_{f_2 \circ f_1}(z) - \sum_{f_2(y)=z}\frac{1}{|Y|} \right|
    &\leq \sum_{z \in Z} \sum_{f_2(y)=z} \left| P_{f_1}(y) - \frac{1}{|Y|}
    \right| \\
    &= \sum_{y \in Y} \left| P_{f_1}(y) - \frac{1}{|Y|} \right| \leq
    \delta_1,
\end{align*}
and by the assumption on $f_2$,
\[
    \sum_{z \in Z} \left|\sum_{f_2(y)=z}\frac{1}{|Y|} - \frac{1}{|Z|}
    \right| = \sum_{z \in Z} \left| \frac{|f_2^{-1}(z)|}{|Y|} -
    \frac{1}{|Z|} \right| = \sum_{z \in Z} \left|
    P_{f_2}(z) - \frac{1}{|Z|} \right| \leq \delta_2.
\]

Therefore,
\begin{align*}
    \left\| P_{f_2 \circ f_1} - U_Z \right\|_1 &=
    \sum_{z \in Z} \left|
    P_{f_2 \circ f_1}(z) - \frac{1}{|Z|} \right| \\
    &\leq \sum_{z \in Z} \left(\left|
    P_{f_2 \circ f_1}(z) - \sum_{f_2(y)=z}\frac{1}{|Y|} \right| +
    \left| \sum_{f_2(y)=z}\frac{1}{|Y|} - \frac{1}{|Z|} \right|\right) \\
    & \leq
    \delta_1 + \delta_2.
\end{align*}
\hal

Let $\d > 0$ and let $f_1,\dots, f_n$ be functions between finite
sets, such that each one of them satisfies $||f_i - U_i||_1 \le \d$,
where $U_i$ is the uniform distribution on the range of $f_i$.
Assume that we apply a finite number of steps $m$, such that each
step is a composition or a direct product of two functions, and
obtain a new function $f$. From the lemmas above we readily deduce
that $f$ satisfies $||f-U||_1 \le m\d$, where $U$ is the uniform
distribution on the range of $f$. By Lemma 3.1(ii), $f$ is therefore
$\sqrt{m\d}$-equidistributed.

Using the above observation and Proposition 1.1 we now easily obtain
the following.

\begin{thm} Let $G$ be a finite group.
Let $m \ge 1$ and $w=[x_1, \dots , x_m]$, an m-fold commutator in
any arrangements of brackets. Then the associated word map $\a_w:G^m
\go G$ is $\g(G)$-equidistributed, where $\g(G) = (m-1)^{1/2}
(\z^G(2)-1)^{1/4}$.

In particular, if $G$ is simple, then $\a_w$ is almost
equidistributed as $|G| \go \infty$.
\end{thm}

Combining Theorems 7.1 and 7.4 with Proposition 3.2 we see that word
maps $\a_w$ associated with $w = x^2y^2$, or with any $m$-fold
commutator $w$, are almost measure preserving on finite simple
groups.

In particular, using the fact that almost all pairs are generating
pairs, we obtain

\begin{cor} Almost all elements $g$ of a finite simple group $G$
can be obtained as $g=[g_1, \dots , g_m]$, an m-fold commutator in
any given arrangements of brackets, where $g_1, \ldots , g_m \in G$
satisfy $\langle g_i, g_j \rangle = G$ for all $i \ne j$.
\end{cor}

We can add various extra conditions on the $m$-tuple $(g_1, \ldots ,
g_m)$ above, provided they hold with probability tending to $1$. For
example, given any non-trivial words $w_1, \ldots , w_k \in F_m$ we
can require that $w_i(g_1, \ldots , g_m) \ne 1$ for all $i = 1,
\ldots ,k$. Indeed, it is proved in \cite{DPSSh} that, if $G$ is a
finite simple group and $1 \ne w \in F_m$, then, as $|G| \go
\infty$, almost all $m$-tuples $(g_1, \ldots , g_m) \in G^m$ satisfy
$w(g_1, \ldots  , g_m) \ne 1$.

For instance, one now easily deduces that almost all elements $g$ of
a finite simple group $G$ can be expressed as $g =
[[g_1,g_2],[g_3,g_4]]$ where  $\langle g_i, g_j \rangle = G$ for $i
\ne j$, and the orders of $g_1, \ldots , g_4$ are as large as we
want.

\end{document}